\documentclass[10pt,twoside,english]{amsart2000}
\advance\oddsidemargin by -1.0cm
\advance\evensidemargin by -1.0cm
\advance\topmargin by -1.0cm 
\usepackage{amssymb}
\usepackage{babel}
\usepackage{amstext2000}
\usepackage{amscd2000}   
\usepackage{epsfig}  
\usepackage{rotating}
\let\mathg\mathfrak
\theoremstyle{plain}
\newtheorem{cor}{Corollary}[section]
\newtheorem{lem}{Lemma}[section]
\newtheorem{thm}{Theorem}[section]

\theoremstyle{definition}
\newtheorem{exa}{Example}[section]
\newtheorem{NB}{Remark}[section]

\newcommand{\bdm}{\begin{displaymath}}
\newcommand{\edm}{\end{displaymath}}
\newcommand{\be}{\begin{equation}}
\newcommand{\ee}{\end{equation}}
\newcommand{\ba}[1]{\begin{array}{#1}}
\newcommand{\ea}{\end{array}}

\newcommand{\btab}{\begin{tabular}}
\newcommand{\etab}{\end{tabular}}

\newcommand{\R}{\ensuremath{\mathbb{R}}}

\newcommand{\End}{\ensuremath{\mathrm{End}}}

\newcommand{\Ad}{\ensuremath{\mathrm{Ad}\,}}

\newcommand{\SU}{\ensuremath{\mathrm{SU}}}

\newcommand{\SO}{\ensuremath{\mathrm{SO}}}
\newcommand{\Spin}{\ensuremath{\mathrm{Spin}}}

\begin{document}
\def\haken{\mathbin{\hbox to 6pt{%
                 \vrule height0.4pt width5pt depth0pt
                 \kern-.4pt
                 \vrule height6pt width0.4pt depth0pt\hss}}}
    \let \hook\intprod
\setcounter{equation}{0}
%
%------ draw title page -----
%
\thispagestyle{empty}
%
%\hbox to \hsize{%
%  \vtop{} \hfill
%  \vtop{\hbox{PRELIMINARY VERSION}}}
%------------------------------
\date{\today}
%----------------------------------------------------------
\title{Spin(9)-structures and connections with totally skew-symmetric torsion}
%----------------------------------------------------------
%
% author and address
%
%-------------------------------------------
%
\author{Thomas Friedrich}
%-------------------------------------------
\address{\hspace{-5mm} 
Thomas Friedrich\newline
Institut f\"ur Reine Mathematik \newline
Humboldt-Universit\"at zu Berlin\newline
Sitz: WBC Adlershof\newline
D-10099 Berlin, Germany\newline
{\normalfont\ttfamily friedric@mathematik.hu-berlin.de}}
%
%------------------------------------------------------
\thanks{Supported by the SFB 288 "Differential geometry
and quantum physics" of the DFG}
%------------------------------------------------------
%----------------------------------------------
\begin{abstract}
%---------------
We study $\Spin(9)$-structures on $16$-dimensional Riemannian manifolds
and characterize the geometric types admitting a connection with totally
skew-symmetric torsion.
\end{abstract}
%-------------
\maketitle
%----------------
\tableofcontents
%----------------
\pagestyle{headings}
%
%
%-------------- body of the document ------------------------------------------
%
%---------------------------------------------------------------------------- 
\section{Introduction}\noindent
%----------------------------------------------------------------------------
The basic model in type II string theory is a $6$-tuple 
$(M^n,g,\nabla,T,\Phi, \Psi)$ consisting of a Riemannian metric $g$,
a metric connection $\nabla$ with totally skew-symmetric torsion form $T$, 
a dilation function $\Phi$ and a spinor field $\Psi$. If the dilation 
function is constant, the string equations can be written in the following 
form (see \cite{Stro} and \cite{IP, Friedrich&I1}):
\bdm
\mbox{Ric}^{\nabla} \ = \ 0, \quad  \delta^g(T) \ = \ 0, 
\quad \nabla \Psi \ = \ 0, \quad T \cdot \Psi \ = \ 0 \, .
\edm
Therefore, an interesting problem is the investigation of metric connections 
with totally skew-symmetric torsion. In \cite{Friedrich&I1} we proved that 
several non-integrable geometric structures (almost contact metric structures,
almost complex structures, $\mbox{G}_2$-structures) admit a unique connection 
$\nabla$ preserving it with totally skew-symmetric torsion. Moreover, we 
computed the corresponding torsion form $T$ and we studied the integrability 
condition for $\nabla$-parallel spinors as well as the Ricci tensor 
$\mbox{Ric}^{\nabla}$. In particular, we constructed $7$-dimensional solutions
of the string equations related to non-integrable $\mbox{G}_2$-structures. The
$5$-dimensional case and its link with contact geometry was investigated 
in more details in the paper \cite{Friedrich&I2}. Similar results concerning 
$8$-dimensional manifolds with a $\Spin(7)$-structure are contained in the 
paper \cite{Iv}, the hyperk\"ahler case was investigated in the papers 
\cite{DoFi}, \cite{IM} and \cite{Ver}. Homogeneous models and the relation 
to Kostant's cubic Dirac operators were discussed in \cite{Agri}. The 
aim of this note is to work out the case of $16$-dimensional Riemannian 
manifolds with a non-integrable $\Spin(9)$-structure. Alfred Gray (see \cite{Gray}) has pointed out that this special geometry may occur
as a geometry with a weak holonomy group. Only recently we once again 
revisited the special $\Spin(9)$-geometries in dimension sixteen and, in 
particular, we proved that there are $4$ basic classes (see \cite{Fri1}). 
Here we will study the problem which of these classes admit a connection 
$\nabla$ with totally skew-symmetric torsion.  
%
%---------------------------------------------------------------------------- 
\section{The geometry of Spin(9)-structures}\noindent
%----------------------------------------------------------------------------
%
The geometric types of $\Spin(9)$-structures on $16$-dimensional oriented 
Riemannian manifolds were investigated in the paper \cite{Fri1}. We
summerize the basic facts defining this special geometry. Let us
consider the $16$-dimensional oriented Euclidean space $\R^{16}$. This 
space is the real spin representation of the group $\Spin(9)$ and, therefore,
there exist nine linear operators $I_{\alpha} : \R^{16} \rightarrow \R^{16}$
such that the following relations hold:
\bdm
I^2_{\alpha} \ = \ \mbox{Id}, \quad I^*_{\alpha} \ = \ I_{\alpha}, \quad 
I_{\alpha} \cdot I_{\beta} + I_{\beta} \cdot I_{\alpha} \ = \ 0 \ \ \
(\alpha \neq \beta), \quad \mbox{Tr}(I_{\alpha}) \ = \ 0 \, .
\edm
The subgroup $\Spin(9) \subset \SO(16)$ can be defined as the group of all
automorphisms of $\R^{16}$ preserving, under conjugation, the $9$-dimensional
subspace $\R^9 := \mbox{Lin}\{I_1, \ldots, I_9\} \subset\mbox{End}(\R^{16})$, 
\bdm
\Spin(9) \ := \ \big\{ g \in \SO(\R^{16}) : \ g \cdot \R^9 \cdot g^{-1} \ = \ 
\R^9 \big\} \, .
\edm
The decomposition of the Lie algebra $\mathg{so}(16) = \mathg{so}(9) \oplus 
\mathg{m}$ is explicitly given by
\bdm
\mathg{so}(9) \ := \ \mbox{Lin}\big\{I_{\alpha} \cdot I_{\beta} \ : \ \alpha < \beta\big\} \ = \ \Lambda^2(\R^9), 
\quad \mathg{m} \ := \ \mbox{Lin}\big\{I_{\alpha} \cdot I_{\beta} \cdot I_{\gamma} \ : \ \alpha < \beta < \gamma\big\} \ = \ \Lambda^3(\R^9) \, .
\edm
The operators $I_{\alpha} \cdot I_{\beta}$ and $I_{\alpha} \cdot I_{\beta} 
\cdot I_{\gamma}$ are skew-symmetric and, consequently, they define two 
systems of $2$-forms $\omega_{\alpha \beta}$ and $\sigma_{\alpha \beta 
\gamma}$.\\

\noindent
Let $(M^{16}, g)$ be an oriented, $16$-dimensional Riemannian manifold. A
$\Spin(9)$-structure is a $9$-dimensional subbundle $V^9 \subset 
\mbox{End}(TM^{16})$ of endomorphisms which is locally generated by sections
$I_{\alpha}$ satisfying the algebraic relations described before. Denote by 
$\mathcal{F}(M^{16})$ the frame bundle of the oriented Riemannian
manifold. Equivalently, a $\Spin(9)$-structure is a reduction $\mathcal{R} 
\subset \mathcal{F}(M^{16})$ of the principal fibre bundle to the subgroup 
$\Spin(9)$. The Levi-Civita connection is a $1$-form on $\mathcal{F}(M^{16})$
with values in the Lie algebra $\mathg{so}(16)$,
\bdm
Z \ : \ T(\mathcal{F}(M^{16})) \longrightarrow \mathg{so}(16) \, .
\edm
We restrict the Levi-Civita connection to a fixed $\Spin(9)$-structure 
$\mathcal{R}$ and decompose it with respect to the decomposition of the
Lie algebra $\mathg{so}(16)$:
\bdm
Z\big|_{T(\mathcal{R})} \ := \ Z^* \, \oplus \ \Gamma \, .
\edm
Then, $Z^*$ is a connection in the principal $\Spin(9)$-bundle $\mathcal{R}$
and $\Gamma$ is a tensorial $1$-form of type Ad, i.\,e., a $1$-form on 
$M^{16}$ with values in the associated bundle
\bdm
\mathcal{R} \times_{\Spin(9)} \mathg{m} \ = \ \mathcal{R} \times_{\Spin(9)} \Lambda^3(\R^9) \ = \ \Lambda^3(V^9) \, .
\edm
The $\Spin(9)$-representation $\R^{16} \otimes \mathg{m} = \R^{16} \otimes 
\Lambda^3(\R^9)$ splits into four irreducible components,
\bdm
\R^{16} \otimes \mathg{m} \ = \ \R^{16} \oplus \mathcal{P}_1(\R^9) \oplus 
\mathcal{P}_2(\R^9) \oplus \mathcal{P}_3(\R^9) \, ,
\edm
and, therefore, we obtain a similar decomposition of the bundle
$\Lambda^1(M^{16}) \otimes \Lambda^3(V^9)$. The representation 
$\mathcal{P}_1(\R^9)$ has dimension $128$. It is the
restriction of the half spin representation $\Delta^{-}_{16}$ of $\Spin(16)$ 
to the subgroup $\Spin(9)$. The dimensions of the irreducible representations
$\mathcal{P}_2(\R^9)$ and $\mathcal{P}_3(\R^9)$ are $432$ and $768$, 
respectively. \\

\noindent
The decomposition of the section $\Gamma$ yields
the classification of all geometric types of $\Spin(9)$-structures. In 
particular, there are four basic classes (see \cite{Fri1}). We remark that 
the sum $\mathcal{P}_1 \oplus \mathcal{P}_2$ is isomorphic to the bundle 
of $3$-forms on $M^{16}$,
\bdm
\Lambda^3(M^{16}) \ = \ \mathcal{P}_1(V^9) \oplus \mathcal{P}_2(V^9) \, .
\edm
In order to fix the normalization, let us describe the embeddings
$\Lambda^i(M^{16}) \longrightarrow \Lambda^1(M^{16}) \otimes \Lambda^3(V^9)$, 
$i = 1,3$, by explicit formulas. If $\mu^1 \in \Lambda^1(M^{16})$ is a 
(co-)vector, then the $1$-form on $M^{16}$ with values in the bundle 
$\Lambda^3(V^9)$ is given by
\bdm
\mu^1 \longmapsto \frac{1}{8}\sum_{\alpha < \beta < 
\gamma}^9 I_{\alpha}I_{\beta}I_{\gamma}(\mu^1) \otimes I_{\alpha} 
\cdot I_{\beta} \cdot I_{\gamma} \, .
\edm
Similarly, if $\mu^3 \in \Lambda^3(M^{16})$ is a $3$-form, we define
\bdm
\mu^3 \longmapsto \frac{1}{8}\sum_{\alpha < \beta < \gamma}^9 
( \sigma_{\alpha\beta\gamma} \haken \mu^3) \otimes I_{\alpha} 
\cdot I_{\beta} \cdot I_{\gamma} \, ,
\edm
where $\sigma_{\alpha\beta\gamma} \haken \mu^3$ denotes the inner product 
of the $2$-forms $\sigma_{\alpha\beta\gamma}$ by $\mu^3$.
%
%---------------------------------------------------------------------------- 
\section{Spin(9)-connections with totally skew-symmetric torsion}\noindent
%----------------------------------------------------------------------------
%
We introduce the following equivariant maps:
\begin{eqnarray*}
\Phi :& {\Bbb R}^{16} \otimes \mathg{spin}(9) \to {\Bbb R}^{16} \otimes S^2 
({\Bbb R}^{16}),& \ \Phi (\Sigma)(X,Y,Z) := g \big(\Sigma (Z)(X),Y \big)+
g \big(\Sigma (Y)(X),Z \big),\\
\Psi :& {\Bbb R}^{16} \otimes \mathg{m} \to {\Bbb R}^{16} \otimes S^2 
({\Bbb R}^{16}),& \
\Psi (\Gamma)(X,Y,Z) := g \big(\Gamma (Y)(X),Z \big)+
g \big(\Gamma (Z)(X),Y \big) \, .
\end{eqnarray*}
It is well known (see \cite{Friedrich&I1}) that a geometric 
$\Spin(9)$-structure admits a connection $\nabla$ with totally skew-symmetric 
torsion if and only if $\Psi(\Gamma)$ is contained in the image of the 
homomorphism $\Phi$. The representation $\R^{16} \otimes \mathg{spin}(9)$ 
splits into
\bdm
\R^{16} \otimes \mathg{spin}(9) \ = \ \R^{16} \oplus \mathcal{P}_1(\R^9) 
\oplus \mathcal{P}_2(\R^9) \, .
\edm
Consequently, if a $\Spin(9)$-structure admits a connection $\nabla$ with 
totally skew-symmetric torsion, then the $\mathcal{P}_3$-part of the form
$\Gamma$ must vanish. We split the $\Spin(9)$-representation 
$\R^{16} \otimes S^2(\R^{16})$ into irreducible components. Since the symmetric
linear maps $I_{\alpha}$ are traceless, the representation $\R^9$ is contained
in $S^2_0(\R^{16})$ and we obtain the decomposition (see \cite{Fri1})
\bdm
\R^{16} \otimes S^2(\R^{16}) \ = \ \R^{16} \oplus \R^{16} \otimes (\R^9 \oplus
D^{126}) \ = \ 2 \cdot \R^{16} \oplus \mathcal{P}_1(\R^9) \oplus \R^{16} 
\otimes  D^{126} \, ,
\edm
where $D^{126}:= \Lambda^4(\R^9)$ is the unique irreducible representation 
of $\Spin(9)$ in dimension $126$. Denote by $D^{672}$ the unique irreducible 
$\Spin(9)$-representation of dimension $672$. Its highest weight is the 
$4$-tuple $(3/2, 3/2, 3/2, 3/2)$.
\begin{lem} 
%----------
The $\Spin(9)$-representation $\R^{16} \otimes S^2(\R^{16})$ splits into 
the irreducible components
\bdm
\R^{16} \otimes S^2(\R^{16}) \ = \ 3 \cdot \R^{16} \oplus 2 \cdot 
\mathcal{P}_1(\R^9) \oplus \mathcal{P}_2(\R^9) \oplus \mathcal{P}_3(\R^9) 
\oplus D^{672} \, .
\edm
\end{lem}
\begin{proof} 
%------------
Since $\R^{16} \otimes \mathg{m}$ contains the representations
$\mathcal{P}_2(\R^9), \ \mathcal{P}_3(\R^9)$ and $\Psi$ is nontrivial, the
tensor product $\R^{16} \otimes D^{126}$ contains the two representations, too.
Moreover, the highest weights of $\R^{16}$ and $D^{126}$ are $(1/2, 1/2, 1/2, 
1/2)$ and $(1, 1, 1, 1)$, respectively. Then the tensor product $\R^{16} 
\otimes D^{126}$ contains the representation $D^{672}$ of highest weight
$(3/2, 3/2, 3/2, 3/2)$ (see \cite{FH}, page 425). Consequently, we obtain
\bdm
\R^{16} \otimes D^{126} \ = \ \mathcal{P}_2(\R^9) \oplus  \mathcal{P}_3(\R^9) 
\oplus D^{672} \oplus S \, ,
\edm
where the dimension of the rest equals $\mbox{dim}(S) = 144$. The
representation $S$ is not an $\SO(9)$-representation. The list of 
small-dimensional $\Spin(9)$-representations yields that 
$S = \R^{16} \oplus  \mathcal{P}_1(\R^9)$, the final result. The 
decomposition of $\R^{16} \otimes D^{126}$ can be computed by a suitable
computer program, too.
\end{proof} 
\begin{lem}\label{Lem1} 
%---------------------
For any two vectors $X, Y \in \R^{16}$ the following identity holds:
\bdm
\sum_{\alpha < \beta}^9 \omega_{\alpha\beta}(X,Y) \cdot \omega_{\alpha\beta} + 
\sum_{\alpha < \beta < \gamma}^9 \sigma_{\alpha\beta\gamma}(X,Y) \cdot \sigma_{\alpha\beta\gamma} \ = \ 8 \cdot X \wedge Y \, .
\edm
\end{lem}
\begin{proof} 
%------------
The $2$-forms $\omega_{\alpha\beta}$ and $\sigma_{\alpha\beta\gamma}$ constitute a basis of the space $\Lambda^2(\R^{16})$ of all $2$-forms in 
sixteen variables. Therefore, the identity is simply the decomposition of the
$2$-form $X \wedge Y$ with respect to this basis. Remark that the length of
the basic forms $\omega_{\alpha\beta}$ and $\sigma_{\alpha\beta\gamma}$ equals
$2 \cdot \sqrt{2}$.
\end{proof}
\begin{thm}\label{thm1} 
%----------------------
A $\Spin(9)$-structure on a $16$-dimensional Riemannian manifold $M^{16}$ 
admits a connection $\nabla$ with totally skew-symmetric torsion if and only 
if the $(\R^{16} \oplus \mathcal{P}_3)$-part of the form $\Gamma$ vanishes. In this case $\Gamma$ is a usual $3$-form on the manifold $M^{16}$, the 
connection $\nabla$ is unique and its torsion form $T$ is given by the formula
$T = - \, 2 \cdot \Gamma$.
\end{thm} 
\begin{proof} 
%------------
For a fixed vector $\Gamma \in \R^{16}$ the tensor 
$\Psi(\Gamma)(X,Y,Y)$ is given by the formula
\bdm
\Psi(\Gamma)(X,Y,Y) \ = \ \frac{1}{4}  \sum_{\alpha < \beta < \gamma}^9 \sigma_{\alpha\beta\gamma}(\Gamma,Y) \cdot \sigma_{\alpha\beta\gamma}(X,Y) \, .
\edm
Since the multiplicity of $\R^{16}$ in the representation $\R^{16} \otimes 
\mathg{spin}(9)$ equals one, any $\Spin(9)$-equivariant map $\Sigma : 
\R^{16} \rightarrow \R^{16} \otimes \mathg{spin}(9)$ is a multiple of
\bdm
\Sigma(\Gamma) \ = \ \sum_{\alpha < \beta}^9 I_{\alpha\beta}(\Gamma) \otimes
I_{\alpha\beta} \, .
\edm
Consequently, if $\Psi(\Gamma)$ is in the image of $\Phi$, there 
exists a constant $c$ such that
\bdm
\sum_{\alpha < \beta < \gamma}^9\sigma_{\alpha\beta\gamma}(\Gamma,Y) \cdot 
\sigma_{\alpha\beta\gamma}(X,Y) \ = \ c \cdot \sum_{\alpha < 
\beta}^9\omega_{\alpha\beta}(\Gamma,Y) \cdot \omega_{\alpha\beta}(X,Y) \, .
\edm 
For $\Gamma = X = e_{16}$ we compute the corresponding
quadratic forms in the variables $y_1, \ldots , y_{16}$:
\bdm
\Psi(e_{16}) \ = \ \sum_{i=1}^8 y_i^2 \, + \, 4 \cdot 
\sum_{j=9}^{15} y_j^2, \quad \quad 
\Phi(\Sigma(e_{16})) \ = \ 7 \cdot \sum_{i=1}^8 y_i^2 \, + \, 
4 \cdot \sum_{j=9}^{15} y_j^2 \, ,
\edm
a contradiction. Next consider the case that $\Gamma \in \Lambda^3(\R^{16})$ is
a $3$-form. By Lemma \ref{Lem1}  we have
\bdm
\Psi(\Gamma)(X,Y,Y) \ = \ \frac{1}{4} \sum_{\alpha\beta\gamma}^9
\Gamma(\sigma_{\alpha\beta\gamma}, Y) \cdot \sigma_{\alpha\beta\gamma}(X, Y) 
\ = \ 
- \, \frac{1}{4} \sum_{\alpha\beta}^9 \Gamma(\omega_{\alpha\beta}, Y) 
\cdot \omega_{\alpha\beta}(X, Y) + 2 \cdot \Gamma(X, Y, Y) \, .
\edm
Since $\Gamma$ is a $3$-form, the term $\Gamma(X,Y,Y)$ vanishes. Let us 
introduce
\bdm
\Sigma(\Gamma) \ := \ - \, \frac{1}{8}\sum_{\alpha\beta}^9 
(\omega_{\alpha\beta} \haken\Gamma) \otimes \omega_{\alpha\beta} \, .
\edm
Then $\Sigma(\Gamma)$ belongs to the space $\R^{16} \otimes \mathg{spin}(9)$ 
and we have $\Phi(\Sigma(\Gamma)) \ = \ \Psi(\Gamma)$. Consequently, in case 
$\Gamma$ is a $3$-form on $M^{16}$, there exists a unique connection 
$\nabla$ preserving the $\Spin(9)$-structure with totally skew-symmetric 
torsion. Its torsion form $T$ is basically given by the difference 
$\Gamma(X) - \Sigma(\Gamma)(X)$ (see \cite{Friedrich&I1}) and we obtain the 
formula $T = - \, 2 \cdot \Gamma$.
\end{proof}

\noindent
Let us characterize $\Spin(9)$-structures of type $\mathcal{P}_1 \oplus 
\mathcal{P}_2$ using the Riemannian covariant derivatives $\nabla I_{\alpha}$ 
of the symmetric endomorphisms describing the structure. For an arbitrary 
$2$-form $S$ we introduce the symmetric forms by the formula
\bdm
S_{\alpha}(Y,Z) \ := \ - \, S(I_{\alpha}(Y), Z) \, + \, S(Y, I_{\alpha}(Z)), \quad \alpha = 1, \ldots 9 \, .
\edm
The connection $\nabla$ preserves the $9$-dimensional bundle of endomorphisms
$I_{\alpha}$ and therefore there exist $1$-forms $M_{\alpha\beta}$ such that
$\nabla I_{\alpha} \ = \ \sum_{\beta=1}^9 M_{\alpha\beta} \cdot I_{\beta}$.
Since $\nabla_XY = \nabla_X^g + \frac{1}{2} \cdot T(X,Y,.)$ we obtain the 
following formula for the Riemannian covariant derivative of the endomorphisms
 $I_{\alpha}$
\bdm
\nabla_X^g I_{\alpha} \ = \ \sum_{\beta=1}^9 M_{\alpha\beta}(X) 
\cdot I_{\beta} + \frac{1}{2} \cdot (X \haken T)_{\alpha} \, ,
\edm
where $T$ is a $3$-form. The latter equation characterizes 
$\Spin(9)$-structures of type $\mathcal{P}_1 \oplus \mathcal{P}_2$.
%
%---------------------------------------------------------------------------- 
\section{Homogeneous Spin(9)-structures}\noindent
%----------------------------------------------------------------------------
%
Consider a Lie group $G$, a subgroup $H$ and suppose that the homogeneous space
$G/H$ is naturally reductive of dimension $16$. We fix a decomposition
\bdm
\mathg{g} \ = \ \mathg{h} \oplus \mathg{n}, \quad [\mathg{h},\mathg{n}] \subset
\mathg{n}, \quad \mathg{n} \ = \ \R^{16}
\edm
as well as a scalar product $(\ , \ )_{\mathg{n}}$ such that for all 
$X,Y,Z \in \mathg{n}$
\bdm
\big([X,Y]_{\mathg{n}}\, , \, Z\big)_{\mathg{n}}+ \big(Y\, , \, [X,Z]_{\mathg{n}}\big)_{\mathg{n}} \ = \ 0
\edm
holds, where $[X,Y]_{\mathg{n}}$ denotes the $\mathg{n}$-part of the 
commutator. Moreover, suppose that the isotropy representation leaves a 
$\Spin(9)$-structure in the vector space $\mathg{n}$ invariant. Then $G/H$ 
admits a homogeneous $\Spin(9)$-structure. Indeed, the frame bundle is an 
associated bundle,
\bdm
\mathcal{F}(G/H) \ = \ G \times_{\Ad} \SO(\mathg{n}) \ ,
\edm
and $\mathcal{R} := G \rightarrow \mathcal{F}(G/H)$ is a reduction to the 
subgroup $H$ contained in $\Spin(9)$. The canonical connection $\nabla^{can}$ 
of the reductive space preserves the $\Spin(9)$-structure and has totally 
skew-symmetric torsion, 
\bdm
T^{\nabla^{can}}(X,Y,Z) \ = \ - \, \big([X,Y]_{\mathg{n}}\, , \, 
Z \big)_{\mathg{n}} \, .
\edm
Consequently, any homogeneous $\Spin(9)$-structure admits an affine connection
with totally skew-symmetric torsion, i.\,e., it is of type
$\mathcal{P}_1 \oplus \mathcal{P}_2$.
\begin{cor} 
%----------
Any homogeneous $\Spin(9)$-structure on a naturally reductive space
$M^{16} = G/H$ is of type $\mathcal{P}_1 \oplus \mathcal{P}_2$.
\end{cor}
\begin{NB} 
%---------
In particular, for any homogeneous $\Spin(9)$-structure 
the difference $\Gamma$ between the Levi-Civita connection and the
canonical connection is a $3$-form. Indeed, the Levi-Civita connection of a 
reductive space is given by the map $\mathg{n} \rightarrow \End(\mathg{n})$ 
\bdm
X \longmapsto \frac{1}{2} \cdot [ X , \ \cdot]_{\mathg{n}} \, .
\edm
Then we obtain
\bdm
\Gamma(X) \ = \ \frac{1}{2} \cdot \mbox{pr}_{\mathg{m}}
\big([X , \ \cdot \ ]_{\mathg{n}}\big) \ =  \ \frac{1}{32}\sum_{i,j=1}^{16}
\sum_{\alpha<\beta<\gamma}
\big([X,e_i]_{\mathg{n}}, e_j \big)_{\mathg{n}} \cdot 
\sigma_{\alpha\beta\gamma}(e_i,e_j) \cdot \sigma_{\alpha\beta\gamma} \, .
\edm
We write the latter equation in the following form
\bdm
\Gamma(X) \ = \ - \, \frac{1}{16}\sum_{\alpha<\beta<\gamma}
(\sigma_{\alpha\beta\gamma} \haken T^{\nabla^{can}})(X) \cdot \sigma_{\alpha\beta\gamma} \ = \ - \, \frac{1}{2} \cdot T^{\nabla^{can}}(X, \, \cdot \, ,\, \cdot\, ) \, , 
\edm
i.\,e., $\Gamma$ is proportional to the torsion of the canonical connection,
\bdm
\Gamma(X)(Y,Z) \ = \ - \, \frac{1}{2} \cdot T^{\nabla^{can}}(X,Y,Z) \ .
\edm
\end{NB}

\noindent
There are homogeneous $\Spin(9)$-structures on different reductive spaces 
(see \cite{Fri1}).
\begin{exa} 
%----------
The group $\Spin(9)$ acts transitively on 
the sphere $S^{15}$, the isotropy group is isomorphic to $\Spin(7)$ and the 
isotropic representation of the reductive space 
$S^1 \times S^{15} = (S^1 \times \Spin(9))/ \Spin(7)$ is contained in 
$\Spin(9)$. 
\end{exa}
\begin{exa}
The space $S^1 \times S^1 \times 
(\SO(8)/\mbox{G}_2)$ admits a homogeneous $\Spin(9)$-structure. 
\end{exa}
\begin{exa}
The space $\SU(5)/\SU(3)$ admits a homogeneous $\Spin(9)$-structure.
\end{exa}
%
%---------------------------------------------------------------------------- 
\section{$\mbox{G}$-connections with totally skew-symmetric torsion}\noindent
%----------------------------------------------------------------------------
%
The class of $\Spin(9)$-structures corresponding to the representation 
$\R^{16} \subset \R^{16} \otimes \mathg{m}$ is related with conformal changes 
of the metric. Indeed, if $(M^{16}, g, V^9)$ is a Riemannian manifold with a 
fixed $\Spin(9)$-structure $V^9 \subset \mbox{End}(TM^{16})$ and 
$g^* = e^{2 f}\cdot g$ is a conformal change of the metric, then the triple  
$(M^{16}, g^*, V^9)$ is a Riemannian manifold with a $\Spin(9)$-structure, 
too. The fact that the $16$-dimensional class of $\Spin(9)$-structures 
corresponding to $\R^{16}$  is not admissible in Theorem \ref{thm1} means 
that the existence of a connection with totally skew-symmetric torsion and 
preserving a $\Spin(9)$-structure is not invariant under conformal 
transformations of the metric. From this point of view the behavior of 
$\Spin(9)$-structures is different 
from the behavior of $\mbox{G}_2$-structures, $\Spin(7)$-structures, 
quaternionic K\"ahler structures or contact structures 
(see \cite{Friedrich&I2}, \cite{Friedrich&I3}, \cite{Iv}, \cite{IM}). We will 
explain this effect in a more general context.\\

\noindent
Let $\mbox{G} \subset \SO(n)$ be a closed subgroup of the orthogonal group
and decompose the Lie algebra
\bdm
\mathg{so}(n) \ = \ \mathg{g} \oplus \mathg{m} \ . 
\edm
A $\mbox{G}$-structure of a Riemannian manifolds $M^n$ is a reduction 
$\mathcal{R} \subset \mathcal{F}(M^{n})$ of the frame bundle to the subgroup 
$\mbox{G}$. The Levi-Civita connection is a $1$-form $Z$ on 
$\mathcal{F}(M^{n})$ with values in the Lie algebra $\mathg{so}(n)$.
We restrict the Levi-Civita connection to a fixed $\mbox{G}$-structure 
$\mathcal{R}$ and decompose it with respect to the decomposition of the
Lie algebra $\mathg{so}(n)$:
\bdm
Z\big|_{T(\mathcal{R})} \ := \ Z^* \, \oplus \ \Gamma \, .
\edm
Then, $Z^*$ is a connection in the principal $\mbox{G}$-bundle $\mathcal{R}$
and $\Gamma$ is a tensorial $1$-form of type Ad, i.\,e., a $1$-form on $M^{n}$ 
with values in the associated bundle $\mathcal{R} \times_{\mbox{G}} 
\mathg{m}$. The $\mbox{G}$-representation $\R^{n} \otimes \mathg{m}$ splits 
into irreducible components and the corresponding decomposition of $\Gamma$
characterizes the different non-integrable $\mbox{G}$-structures.
We introduce the equivariant maps:
\begin{eqnarray*}
\Phi :& {\Bbb R}^{n} \otimes \mathg{g} \to {\Bbb R}^{n} \otimes S^2 
({\Bbb R}^{n}),& \ \Phi (\Sigma)(X,Y,Z) := g \big(\Sigma (Z)(X),Y \big)+
g \big(\Sigma (Y)(X),Z \big),\\
\Psi :& {\Bbb R}^{n} \otimes \mathg{m} \to {\Bbb R}^{n} \otimes S^2 
({\Bbb R}^{n}),& \
\Psi (\Gamma)(X,Y,Z) := g \big(\Gamma (Y)(X),Z \big)+
g \big(\Gamma (Z)(X),Y \big) \, .
\end{eqnarray*}
It is well known (see \cite{Friedrich&I1}) that a geometric 
$\mbox{G}$-structure admits a connection $\nabla$ with totally skew-symmetric 
torsion if and only if $\Psi(\Gamma)$ is contained in the image of the 
homomorphism $\Phi$. There  is an equivalent formulation of this condition. 
Indeed, let us introduce the maps
\bdm
\Theta_1 : \Lambda^3({\Bbb R}^n) \to {\Bbb R}^{n} \otimes \mathg{m}, \quad
\Theta_2 : \Lambda^3({\Bbb R}^n) \to {\Bbb R}^{n} \otimes \mathg{g}
\edm 
given by the formulas
\bdm
\Theta_1(T) \ := \ \sum_i (\sigma_i \haken T) \otimes \sigma_i, \quad
\Theta_2(T) \ := \ \sum_j (\mu_i \haken T) \otimes \mu_j
\edm
where $\sigma_i$ is an orthonormal basis in $\mathg{m}$ and $\mu_j$ is an
orthonormal basis in $\mathg{g}$. Observe that the kernel of the map 
$(\Psi \oplus \Phi) : {\Bbb R}^{n} \otimes \mathg{so}(n) \to 
{\Bbb R}^{n} \otimes S^2 ({\Bbb R}^{n})$ coincides with the image of the 
map $(\Theta_1 \oplus \Theta_2) : \Lambda^3({\Bbb R}^n) \to {\Bbb R}^{n} 
\otimes \mathg{so}(n)$. Consequently, for any element $\Gamma \in 
{\Bbb R}^{n} \otimes \mathg{m}$, the condition 
$\Psi(\Gamma) \in \mbox{Image}(\Phi)$ is equivalent
to $\Gamma \in \mbox{Image}(\Theta_1)$.
\begin{thm}\label{thm2} 
%----------------------
A $\mbox{G}$-structure $\mathcal{R} \subset \mathcal{F}(M^{n})$ of a 
Riemannian manifold admits a connection $\nabla$ with totally skew-symmetric 
torsion if and only if the $1$-form $\Gamma$ belongs to the image 
of $\Theta_1$, $\Gamma = \Theta_1(T)$. In this case the $3$-form $(- \, 2 
\cdot T)$ is the torsion form of the connection.
\end{thm}
\noindent
Consequently, only such geometric types (i.\,e.\ irreducible components of 
${\Bbb R}^n \otimes \mathg{m}$) are admissible which occur in the 
$\mbox{G}$-decomposition of $\Lambda^3({\Bbb R}^n)$. This explains the different behavior of $\mbox{G}$-structures with respect to conformal transformations.
\begin{exa} 
%----------
In case of $\mbox{G} = \Spin(9)$ we have
\bdm
\R^{16} \otimes \mathg{m} \ = \ \R^{16} \oplus \Lambda^3(\R^{16}) \oplus 
\mathcal{P}_3(\R^9) 
\edm
and the $\R^{16}$-component is \emph{not} contained in $\Lambda^3(\R^{16})$, 
i.\,e., a conformal change of a $\Spin(9)$-structure does \emph{not} 
preserve the property that the structure admits a connection with totally 
skew-symmetric torsion. 
\end{exa}
\begin{exa} 
%----------
In case of a $7$-dimensional $\mbox{G}_2$-structure the situation is different. Indeed, we decompose the $\mbox{G}_2$-representation
(see \cite{Friedrich&I1})
\bdm
\Lambda^3(\R^7) \ = \ \R^1 \oplus \R^7 \oplus \Lambda^3_{27}, \quad 
\R^7 \otimes \mathg{m} \ = \ \R^1 \oplus \R^7 \oplus \Lambda^2_{14} \oplus 
\Lambda^3_{27}  
\edm
and, consequently, a conformal change of a $\mbox{G}_2$-structure preserves 
the property that the structure admits a connection with totally 
skew-symmetric torsion.
\end{exa}
\begin{exa} 
%----------
Let us consider $\Spin(7)$-structures on $8$-dimensional Riemannian
manifolds. The subgroup $\Spin(7) \subset \SO(8)$ is the real 
$\Spin(7)$-representation $\Delta_7 = \R^8$. The complement 
$\mathg{m} = \R^7$ is the standard $7$-dimensional representation and the 
$\Spin(7)$-structures on an $8$-dimensional Riemannian manifold $M^8$ 
correspond to the irreducible components of the tensor product
\bdm
\R^8 \otimes \mathg{m} \ = \ \R^8 \otimes \R^7 \ = \ \Delta_7 \otimes \R^7 \ = \ \Delta_7 \oplus \mbox{K} \, ,
\edm
where $\mbox{K}$ denotes the kernel of the Clifford multiplication 
$\Delta_7 \otimes \R^7 \to \Delta_7$. It is well known that 
$\mbox{K}$ is an irreducible $\Spin$-representation. Therefore, there are 
only two basic types of $\Spin(7)$-structures (see \cite{Fer}). On the other 
hand, the map $\Lambda^3(\R^8) \to \R^8 \otimes 
\mathg{m}$ is injective and the $\Spin(7)$-representation 
$\Lambda^3(\R^8) = \Lambda^3(\Delta_7)$ splits again into the irreducible 
components
\bdm
\Lambda^3(\Delta_7) \ = \ \Delta_7 \oplus \mbox{K} \ ,
\edm
i.\,e.,  $\Lambda^3(\R^8) \to \R^8 \otimes \mathg{m}$ is an isomorphism.
Theorem \ref{thm2} yields immediately that \emph{any $\Spin(7)$-structure on an
$8$-dimensional Riemannian manifold admits a connection with totally skew-symmetric torsion} (see \cite{Iv}). We remark that $n = 8$ is the smallest dimension
where this effect can occur. Indeed, let $\mbox{G} \subset \SO(n)$ be a 
subgroup of dimension $\mbox{g}$ and suppose that any $\mbox{G}$-structure 
admits a connection with totally skew-symmetric torsion, i.e., the map 
 $\Lambda^3(\R^n) \to \R^n \otimes \mathg{m}$ is surjective. On the other side,
 the isotropy representation $\mbox{G} \to \SO(\mathg{m})$ of the compact 
Riemannian manifold $\SO(n)/\mbox{G}$ is injective. Consequently, we obtain
the inequalities
\bdm
\frac{1}{3} ( n^2 - 1) \ \leq \ \mbox{g} \ \leq \ 
\frac{1}{2} (n^2 - 3n +2) \, .
\edm
The minimal pair satisfying this condition is $n = 8, \, \mbox{g} = 21$. Using not only the dimension of the $\mbox{G}$-representation one can exclude
other dimensions, for example $n=9$. 
\end{exa} 
%------------------------------------------
%\addcontentsline{toc}{section}{Literature}
%------------------------------------------
    
\end{document}